\documentclass{article} 

\usepackage{latexsym}

\usepackage{amsmath, amssymb}

\title{A CHARACTERISTIC POLYNOMIAL FOR THE TRANSITION PROBABILITY MATRIX 
OF A CORRELATED RANDOM WALK ON A GRAPH}
\author{Takashi KOMATSU \\
Department of Bioengineering School of Engineering,\\ 
The University of Tokyo \\
Bunkyo, Tokyo, 113-8656, JAPAN \\ 
e-mail: komatsu@coi.t.u-tokyo.ac.jp \\ 
Norio KONNO\\
Department of Applied Mathematics, 
Faculty of Engineering, \\ 
Yokohama National University, \\
Hodogaya, Yokohama 240-8501, JAPAN \\
e-mail: konno@ynu.ac.jp \\
Iwao SATO \\ 
National Institute of Technology, Oyama College, \\ 
Oyama, Tochigi 323-0806, JAPAN \\
e-mail: isato@oyama-ct.ac.jp }
 
\begin{document}
 \maketitle

\clearpage

\begin{abstract}
We define a correlated random walk (CRW) induced from the time evolution matrix 
(the Grover matrix) of the Grover walk on a graph $G$, and present a formula for 
the characteristic polynomial of the transition probability matrix of this CRW 
by using a determinant expression for the generalized weighted zeta function of $G$. 
As applications, we give the spectrum of the transition probability matrices 
for the CRWs induced from the Grover matrices of regular graphs and 
semiregular bipartite graphs. 
Furthermore, we consider another type of the CRW on a graph. 
\end{abstract}

\vspace{5mm}


\vspace{5mm}

The contact author for correspondence:

Iwao Sato 

Oyama National College of Technology, 
Oyama, Tochigi 323-0806, JAPAN 

Tel: +81-285-20-2176

Fax: +81-285-20-2880

E-mail: isato@oyama-ct.ac.jp

\clearpage

\section{Introduction}

Zeta functions of graphs started from the Ihara zeta functions of regular graphs 
by Ihara [6]. 
In [6], he showed that their reciprocals are explicit polynomials. 
A zeta function of a regular graph $G$ associated with a unitary 
representation of the fundamental group of $G$ was developed by 
Sunada [15,16]. 
Hashimoto [4] generalized Ihara's result on the Ihara zeta function of 
a regular graph to an irregular graph, and showed that its reciprocal is 
again a polynomial by a determinant containing the edge matrix. 
Bass [1] presented another determinant expression for the Ihara zeta function 
of an irregular graph by using its adjacency matrix. 

Morita [12] defined a generalized weighted zeta function of a digraph which contains various zeta functions 
of a graph or a digraph. 
Ide et al [5] presented a determinant expression for the above generalized weighted zeta function of a graph.

The time evolution matrix of a discrete-time quantum walk in a graph 
is closely related to the Ihara zeta function of a graph. 
A discrete-time quantum walk is a quantum analog of the classical random walk on a graph whose state vector is governed by 
a matrix called the time evolution matrix(see [8]). 
Ren et al. [13] gave a relationship between the discrete-time quantum walk and the Ihara zeta function of a graph.  
Konno and Sato [10] obtained a formula of the characteristic polynomial of the Grover matrix 
by using the determinant expression for the second weighted zeta function of a graph.

In this paper, we define introduce a new correlated random walk induced from the time evolution matrix (the Grover matrix) 
of the Grover walk on a graph, and present a formula for the characteristic polynomial of its transition probability matrix.

In Section 2, we review for the Ihara zeta function and the generalized weighted zeta functions of a graph. 
In Section 3, we review for the Grover walk on a graph. 
In Section 4, we define a correlated random walk (CRW) induced from the time evolution matrix (the Grover matrix) of 
the Grover walk on a graph $G$, and present a formula for the characteristic polynomial of the transition probability matrix 
of this CRW. 
In Section 5, we give the spectrum of the transition probability matrix for this CRW of a regular graph. 
In Section 6, we present the spectrum for the transition probability matrix of 
this CRW of a semiregular bipartite graph. 
In Section 7, we present formulas for the characteristic polynomials of the transition probability matrices  
of another type of the CRW on a graph, and give the spectrum of its transition probability matrix.

\section{Preliminaries} 

\subsection{Zeta functions of graphs}

Graphs and digraphs treated here are finite.
Let $G$ be a connected graph and $D_G$ the symmetric digraph 
corresponding to $G$. 
Set $D(G)= \{ (u,v),(v,u) \mid uv \in E(G) \} $. 
For $e=(u,v) \in D(G)$, set $u=o(e)$ and $v=t(e)$. 
Furthermore, let $e^{-1}=(v,u)$ be the {\em inverse} of $e=(u,v)$. 
For $v \in V(G)$, the {\em degree} $\deg {}_G v = \deg v = d_v $ is the number of vertices 
adjacent to $v$ in $G$.  
A graph $G$ is called {\em $k$-regular} if $\deg v=k$ for each $v \in V(G)$. 

A {\em path $P$ of length $n$} in $G$ is a sequence 
$P=(e_1, \cdots ,e_n )$ of $n$ arcs such that $e_i \in D(G)$,
$t( e_i )=o( e_{i+1} )(1 \leq i \leq n-1)$. 
If $e_i =( v_{i-1} , v_i )$ for $i=1, \cdots , n$, then we write 
$P=(v_0, v_1, \cdots ,v_{n-1}, v_n )$. 
Set $ \mid P \mid =n$, $o(P)=o( e_1 )$ and $t(P)=t( e_n )$. 
Also, $P$ is called an {\em $(o(P),t(P))$-path}. 
We say that a path $P=( e_1 , \cdots , e_n )$ has a {\em backtracking} 
if $ e^{-1}_{i+1} =e_i $ for some $i\ (1 \leq i \leq n-1)$. 
A $(v, w)$-path is called a {\em $v$-cycle} 
(or {\em $v$-closed path}) if $v=w$. 
The {\em inverse cycle} of a cycle 
$C=( e_1, \cdots ,e_n )$ is the cycle 
$C^{-1} =( e^{-1}_n, \cdots ,e^{-1}_1 )$.

We introduce an equivalence relation between cycles. 
Two cycles $C_1 =(e_1, \cdots ,e_m )$ and 
$C_2 =(f_1, \cdots ,f_m )$ are called {\em equivalent} if there exists a positive number $k$ 
such that $f_j =e_{j+k} $ for all $j$, where the subscripts are considered by modulo $m$. 
The inverse cycle of $C$ is in general not equivalent to $C$. 
Let $[C]$ be the equivalence class which contains a cycle $C$. 
Let $B^r$ be the cycle obtained by going $r$ times around a cycle $B$. 
Such a cycle is called a {\em multiple} of $B$. 
A cycle $C$ is {\em reduced} if 
both $C$ and $C^2 $ have no backtracking. 
Furthermore, a cycle $C$ is {\em prime} if it is not a multiple of 
a strictly smaller cycle. 
Note that each equivalence class of prime, reduced cycles of a graph $G$ 
corresponds to a unique conjugacy class of 
the fundamental group $ \pi {}_1 (G,v)$ of $G$ at a vertex $v$ of $G$. 

The {\em Ihara(-Selberg) zeta function} of $G$ is defined by  
\[
{\bf Z} (G,u)= \prod_{[C]} (1- u^{ \mid C \mid } )^{-1} , 
\]
where $[C]$ runs over all equivalence classes of prime, reduced cycles 
of $G$. 

Let $G$ be a connected graph with $n$ vertices and $m$ edges. 
Then two $2m \times 2m$ matrices 
${\bf B} = {\bf B} (G)=( {\bf B}_{e,f} )_{e,f \in D(G)} $ and 
${\bf J}_0 ={\bf J}_0 (G) =( {\bf J}_{e,f} )_{e,f \in D(G)} $ 
are defined as follows: 
\[
{\bf B}_{e,f} =\left\{
\begin{array}{ll}
1 & \mbox{if $t(e)=o(f)$, } \\
0 & \mbox{otherwise}
\end{array}
\right.
, 
{\bf J}_{e,f} =\left\{
\begin{array}{ll}
1 & \mbox{if $f= e^{-1} $, } \\
0 & \mbox{otherwise.}
\end{array}
\right.
\]
The matrix ${\bf B} - {\bf J}_0 $ is called the {\em edge matrix} of $G$. 

\newtheorem{theorem}{Theorem}
\begin{theorem}[Hashimoto; Bass] 
Let $G$ be a connected graph with $n$ vertices and $m$ edges. 
Then the reciprocal of the Ihara zeta function of $G$ is given by 
\[
{\bf Z} (G,u) {}^{-1} = \det ( {\bf I}_{2m} - u ( {\bf B} - {\bf J}_0 )) 
=(1- u^2 )^{m-n} \det( {\bf I}_n  -u {\bf A} (G)+ u^2 ( {\bf D}_G - {\bf I}_n  )) , 
\]
where ${\bf D}_G =( d_{ij} )$ is the diagonal matrix with 
$d_{ii} = \deg {}_G \  v_i \  (V(G)= \{ v_1 , \cdots , v_n \} )$. 
\end{theorem}

The first identity in Theorem 1 was obtained by Hashimoto [4]. 
Also, Bass [1] proved the second identity by using a linear algebraic method. 

Stark and Terras [14] gave an elementary proof of this formula, and 
discussed three different zeta functions of any graph. 
Various proofs of Bass' Theorem were given by Kotani and Sunada [11], 
and Foata and Zeilberger [3].

\subsection{The generalized weighted zeta functions of a graph}

Let $G$ be a connected graph with $n$ vertices and $m$ edges, and 
$D(G)= \{ e_1, \ldots, e_m, e_{m+1} , \ldots , $ 
$e_{2m} \} (e_{m+i} = e^{-1}_i (1 \leq i \leq m ))$. 
Furthermore, we consider two functions $\tau : D(G) \longrightarrow \mathbb{C}$ and 
$\mu : D(G) \longrightarrow \mathbb{C}$. 
Let $\theta : D(G) \times D(G) \longrightarrow \mathbb{C} $ be a function such that 
\[
\theta (e,f)= \tau (f) \delta {}_{t(e)o(f)} - \mu (f) \delta {}_{e^{-1} f} .  
\]
We introduce a $2m \times 2m$ matrix ${\bf M} ( \theta ) =( M_{ef} )_{e,f \in D(G)} $ as follows: 
\[
M_{ef} =\theta (e,f) . 
\]
Then the {\em generalized weighted zeta function} ${\bf Z} {}_G (u, \theta )$ of $G$ is defined 
as follows(see [12]):  
\[
{\bf Z} {}_G (u, \theta )= \det ( {\bf I}_{2m} -u {\bf M} ( \theta ) )^{-1} .  
\]

We consider two $n \times n$ matrices ${\bf A}_G ( \theta )=( a_{uv} )_{u, v \in V(G)} $ and 
${\bf D}_G ( \theta )=( d_{uv} )_{u, v \in V(G)} $ as follows: 
\[
a_{uv} =\left\{
\begin{array}{ll}
\tau (e)/(1- u^2 \mu (e) \mu (e^{-1} )) & \mbox{if $e(u,v) \in D(G)$, } \\
0 & \mbox{otherwise, }
\end{array}
\right. 
\] 
\[
d_{uv} =\left\{
\begin{array}{ll}
\sum_{o(e)=u} \tau (e) \mu ( e^{-1} )/(1- u^2 \mu (e) \mu (e^{-1} )) & \mbox{if $u=v$, } \\
0 & \mbox{otherwise.}
\end{array}
\right. 
\]

A determinant expression for the generalized weighted zeta function of a graph is given as follows(see [5]):

\begin{theorem}[Ide, Ishikawa, Morita, Sato and Segawa]
Let $G$ be a connected graph with $n$ vertices and $m$ edges, and let $\tau : D(G) \longrightarrow \mathbb{C}$ and 
$\mu : D(G) \longrightarrow \mathbb{C}$ be two functions.  
Then 
\[
{\bf Z} {}_G (u, \theta ) {}^{-1} = \prod^m_{j=1} (1- u^2 \mu (e_j ) \mu (e^{-1}_j )) 
\det ( {\bf I}_{n} -u {\bf A}_G ( \theta )+ u^2 {\bf D}_G ( \theta )) ,  
\]
where $D(G)= \{ e_1 , \ldots , e_m , e_{m+1} , \ldots , e_{2m} \} \ (e_{m+j} = e^{-1}_j (1 \leq j \leq m))$. 
\end{theorem}

\section{The Grover walk on a graph}

Let $G$ be a connected graph with $n$ vertices and $m$ edges, 
$V(G)= \{ v_1 , \ldots , v_n \} $ and $D(G)= \{ e_1 , \ldots , e_m , 
e^{-1}_1 , \ldots , e^{-1}_m \} $. 
Set $d_j = d_{v_j} = \deg v_j $ for $i=1, \ldots ,n$. 
The {\em Grover matrix} ${\bf U} ={\bf U} (G)=( U_{ef} )_{e,f \in R(G)} $ 
of $G$ is defined by 
\[
U_{ef} =\left\{
\begin{array}{ll}
2/d_{t(f)} (=2/d_{o(e)} ) & \mbox{if $t(f)=o(e)$ and $f \neq e^{-1} $, } \\
2/d_{t(f)} -1 & \mbox{if $f= e^{-1} $, } \\
0 & \mbox{otherwise.}
\end{array}
\right. 
\]
The discrete-time quantum walk with the matrix ${\bf U} $ as a time evolution matrix 
is called the {\em Grover walk} on $G$.

Let $G$ be a connected graph with $n$ vertices and $m$ edges. 
Then the $n \times n$ matrix ${\bf T} (G)=( T_{uv} )_{u,v \in V(G)}$ is given as follows: 
\[
T_{uv} =\left\{
\begin{array}{ll}
1/( \deg {}_G u)  & \mbox{if $(u,v) \in D(G)$, } \\
0 & \mbox{otherwise.}
\end{array}
\right.
\] 
Note that the matrix ${\bf T} (G)$ is the transition matrix of the simple random walk on $G$(see [10]).

\begin{theorem}[Konno and Sato]
Let $G$ be a connected graph with $n$ vertices $v_1 , \ldots , v_n $ and $m$ edges. 
Then the characteristic polynomial for the Grover matrix ${\bf U}$ of $G$ is given by 
\[
\begin{array}{rcl}
\det ( \lambda {\bf I}_{2m} - {\bf U} ) & = & ( \lambda {}^2 -1)^{m-n} \det (( \lambda {}^2 +1) {\bf I}_n -2 \lambda {\bf T} (G)) \\
\  &   &                \\ 
\  & = & \frac{( \lambda {}^2 -1)^{m-n} \det (( \lambda {}^2 +1) {\bf D} -2 \lambda {\bf A} (G))}{d_{v_1} \cdots d_{v_n }} .   
\end{array}
\]
\end{theorem}

From this Theorem, the spectra of the Grover matrix on a graph is obtained by means of those of 
${\bf T} (G)$ (see [13]). 
Let $Spec ({\bf F})$ be the spectra of a square matrix ${\bf F}$ .

\newtheorem{corollary}{Corollary} 
\begin{corollary}[Emms, Hancock, Severini and Wilson] 
Let $G$ be a connected graph with $n$ vertices and $m$ edges. 
The Grover matrix ${\bf U}$ has $2n$ eigenvalues of the form 
\[
\lambda = \lambda {}_T \pm i \sqrt{1- \lambda {}^2_T } , 
\]
where $\lambda {}_T $ is an eigenvalue of the matrix ${\bf T} (G)$. 
The remaining $2(m-n)$ eigenvalues of ${\bf U}$ are $\pm 1$ with equal multiplicities. 
\end{corollary}

\section{A correlated random walk on a graph}

Let $G$ be a connected graph with $n$ vertices and $m$ edges, and 
${\bf U} $ be the Grover matrix of $G$.  
Then we define a $2m \times 2m$ matrix ${\bf P} =( P_{ef} )_{e, f \in D(G)} $ as follows: 
\[
P_{ef} =| U_{ef} |^2 . 
\]
Note that 
\[ 
P_{ef} = \left\{
\begin{array}{ll}
4/d^2_{t(f)} (=4/d^2_{o(e)} ) & \mbox{if $t(f)=o(e)$ and $f \neq e^{-1} $, } \\
( 2/d_{t(f)} -1)^2 & \mbox{if $f= e^{-1} $, } \\
0 & \mbox{otherwise.}
\end{array}
\right. 
\]
The random walk with the matrix ${\bf P} $ as a transition probability matrix is called the {\em correlated random walk (CRM)} 
(with respect to the Grover matrix) on $G$(see [7,9]). 

Let ${\bf R} =( R_{ef} )_{e,f \in D(G)}$ be a $2m \times 2m$ matrix such that 
\[ 
R_{ef} = \left\{
\begin{array}{ll}
4/d^2_{o(f)} (=4/d^2_{o(e)} ) & \mbox{if $o(e)=o(f)$ and $f \neq e$, } \\
( 2/d_{o(f)} -1)^2 & \mbox{if $f=e$, } \\
0 & \mbox{otherwise.}
\end{array}
\right. 
\]
Then we have 
\[
{\bf P} = {\bf J}_0 {\bf R} . 
\]

By Theorem 2, we obtain the following formula for ${\bf P} $.

\begin{theorem}
Let $G$ be a connected graph with $n$ vertices and $m$ edges,  
and let ${\bf P} $ be the transition probability matrix of the CRW with respect to the Grover matrix. 
Then 
\[
\det ( {\bf I}_{2m} -u {\bf P} )= \prod^m_{j=1} (1- u^2 (\frac{4}{d_{o(e_j )}} -1)( \frac{4}{d_{t(e_j)}} -1)) 
\det ( {\bf I}_n -u {\bf A}_{CRW} + u^2 {\bf D}_{CRW} ) , 
\]
where 
\[
( {\bf A}_{CRW} )_{xy}= 
\left\{
\begin{array}{ll} 
\frac{4/d^2_x }{1- u^2 (4/d_x -1)(4/d_y -1)} & \mbox{if $(x,y) \in D(G)$, } \\
0 & \mbox{otherwise, }  
\end{array}
\right. 
\]
\[ 
( {\bf D}_{CRW} )_{xy}= 
\left\{
\begin{array}{ll} 
\sum_{o(e)=x} \frac{4/d^2_x (4/d_{t(e)} -1)}{1- u^2 (4/d_x -1)(4/d_{t(e)} -1)} & \mbox{if $x=y$, } \\
0 & \mbox{otherwise.}
\end{array}
\right. 
\] 
\end{theorem}

{\bf Proof}.  For the matrix ${\bf P}$, we have 
\[
P_{ef} = \frac{4}{ d^2_{o(e)}} \delta {}_{t(f)o(e)} - ( \frac{4}{d_{o(e)}} -1) \delta {}_{f^{-1} e} . 
\]
The we let two functions $\tau : D(G) \longrightarrow \mathbb{C}$ and $\mu : D(G) \longrightarrow \mathbb{C}$.  
as follows: 
\[ 
\tau (e)= \frac{4}{d^2_{o(e)}} \ and \ \mu (e)= \frac{4}{d_{o(e)}} -1 . 
\]
Furthermore, let 
\[
\theta (e,f)= \frac{4}{d^2_{o(f)}} \delta {}_{t(e)o(f)} -( \frac{4}{d_{o(f)}} -1) \delta {}_{e^{-1} f}  . 
\]
Then we have 
\[ 
{\bf P} = {}^t {\bf M} ( \theta ) . 
\]
Thus, we obtain 
\[
\det ( {\bf I}_{2m} -u {\bf P} )= \det ( {\bf I}_{2m} -u \ {}^t {\bf M} ( \theta ))= \det ( {\bf I}_{2m} -u {\bf M} ( \theta )) 
= {\bf Z} {}_G (u, \theta )^{-1} .   
\]
By Theorem  , we have 
\[
\det ( {\bf I}_{2m} -u {\bf P} )= \prod^m_{j=1} (1- u^2 (\frac{4}{d_{o(e_j )}} -1)( \frac{4}{d_{t(e_j)}} -1)) 
\det ( {\bf I}_n -u {\bf A}_{CRW} + u^2 {\bf D}_{CRW} ) , 
\]
where 
\[
( {\bf A}_{CRW} )_{xy}= 
\left\{
\begin{array}{ll} 
\frac{4/d^2_x }{1- u^2 (4/d_x -1)(4/d_y -1)} \mbox{if $(x,y) \in D(G)$, } \\
0 & \mbox{otherwise, }  
\end{array}
\right. 
\]
\[ 
( {\bf D}_{CRW} )_{xy}= 
\left\{
\begin{array}{ll} 
\sum_{o(e)=x} \frac{4/d^2_x (4/d_{t(e)} -1)}{1- u^2 (4/d_x -1)(4/d_{t(e)} -1)} & \mbox{if $x=y$, } \\
0 & \mbox{otherwise.}
\end{array}
\right. 
\] 
$\Box$

By Theorem 4, we obtain the spectrum of the transition probability matrices  
for the CRWs induced from the Grover matrices of regular graphs and semiregular bipartite graphs.

\section{An application to the correlated random walk on a regular graph} 

We present spectra for the transition matrix of the correlated random walk on a regular graph with respect to 
the Grover matrix.

\begin{theorem}
Let $G$ be a connected $d$-regular graph with $n$ vertices and $m$ edges, where $d \geq 2$. 
Furthermore, let ${\bf P} $ be the transition probability matrix of the CRW with respect to the Grover matrix. 
Then 
\[
\det ({\bf I}_{2m} -u {\bf P} )= \frac{(d^2 - u^2 (4-d)^2 )^{m-n} }{d^{2m}} \det (d(d+(4-d) u^2 ) {\bf I}_n -4u {\bf A} (G)) . 
\]
\end{theorem}

{\bf Proof}.  Let $G$ be a connected $d$-regular graph with $n$ vertices and $m$ edges, where $d \geq 2$. 
Then we have 
\[
d_{o(e)} = d_{t(e)} =d \ for \ each \ e \in D(G) . 
\]
Thus, we have 
\[
1-u^2 ( \frac{4}{ d_{o(e)}} -1)( \frac{4}{ d_{t(e)}} -1)= \frac{d^2 - u^2 (4-d)^2}{d^2 } , 
\]
\[
( {\bf A}_{CRW} )_{xy} = \frac{4/d^2_x }{1- u^2 (4/d_x -1)(4/d_y -1)}= \frac{4}{d^2 - u^2 (4-d)^2} \ if \ (x,y) \in D(G) 
\]
and 
\[ 
( {\bf D}_{CRW} )_{xy}= \sum_{o(e)=x} \frac{4/d^2_x (4/d_{t(e)} -1)}{1- u^2 (4/d_x -1)(4/d_{t(e)} -1)} 
=d \cdot \frac{4(4-d)}{d(d^2 - u^2 (4-d)^2 )} = \frac{4(4-d)}{d^2 - u^2 (4-d)^2} . 
\]
Therefore, it follows that 
\[
{\bf A}_{CRW} = \frac{4}{d^2 - u^2 (4-d)^2} {\bf A} (G) \ and \ {\bf D}_{CRW} = \frac{4(4-d)}{d^2 - u^2 (4-d)^2} {\bf I}_n . 
\]

By Theorem 4, we have 
\[
\begin{array}{rcl} 
\  &   & \det ({\bf I}_{2m} -u {\bf P} ) \\
\  &   &                \\ 
\  & = & \frac{(d^2 - u^2 (4-d)^2 )^m }{d^{2m}} \det ( {\bf I}_n -u \frac{4}{d^2 - u^2 (4-d)^2} {\bf A} (G)
+ u^2 \frac{4(4-d)}{d^2 - u^2 (4-d)^2} {\bf I}_n ) \\
\  &   &                \\ 
\  & = & \frac{(d^2 - u^2 (4-d)^2 )^{m-n} }{d^{2m}} \det ((d^2 - u^2 (4-d)^2 ) {\bf I}_n -4u {\bf A} (G)+4(4-d) u^2 {\bf I}_n ) \\
\  &   &                \\ 
\  & = & \frac{(d^2 - u^2 (4-d)^2 )^{m-n} }{d^{2m}} \det (d(d+(4-d) u^2 ) {\bf I}_n -4u {\bf A} (G)) . 
\end{array}
\] 
$\Box$

By substituting $u=1/ \lambda $, we obtain the following result.

\begin{corollary}
Let $G$ be a connected $d$-regular graph with $n$ vertices and $m$ edges, where $d \geq 2$. 
Furthermore, let ${\bf P} $ be the transition probability matrix of the CRW with respect to the Grover matrix. 
Then 
\[
\begin{array}{rcl}
\det ( \lambda {\bf I}_{2m} - {\bf P} ) & = & \frac{(d^2 \lambda {}^2 -(4-d)^2 )^{m-n} }{d^{2m}} \det (d(d \lambda {}^2 
+(4-d)) {\bf I}_n -4 \lambda {\bf A} (G)) \\  
\  &   &                \\ 
\  & = & ( \lambda {}^2 -( \frac{4}{d} -1)^2 )^{m-n} \lambda {}^n \det (( \lambda +( \frac{4}{d} -1) ) \frac{1}{ \lambda } ) 
{\bf I}_n - \frac{4}{ d^2 } {\bf A} (G)) . 
\end{array}
\]
\end{corollary}

The second identity of Corollary 2 is considered as the spectral mapping theorem for ${\bf P} $. 

By Corollary 2, we obtain the spectra for the transition matrix ${\bf P} $ of the CRW 
with respect to the Grover matrix on a regular graph.

\begin{corollary}
Let $G$ be a connected $d( \geq 2)$-regular graph with $n$ vertices and $m$ edges. 
Then the transition probability matrix ${\bf P} $ has $2n$ eigenvalues of the form 
\[ 
\lambda = \frac{2 \lambda {}_A \pm \sqrt{4 \lambda {}_A {}^2 - d^3 (4-d)}}{ d^2 } , 
\]
where $\lambda {}_A $ is an eigenvalue of the matrix ${\bf A} (G)$. 
The remaining $2(m-n)$ eigenvalues of ${\bf P} $ are $\pm (4-d)/d$ with equal multiplicities $m-n$. 
\end{corollary}

{\bf Proof}. By Corollary 2, we have 
\[
\begin{array}{rcl} 
\det ( \lambda {\bf I}_{2m} - {\bf P} ) & = & ( d^2  \lambda {}^2 -(4-d)^2 )^{m-n} / d^{2m} 
\prod_{ \lambda {}_A \in Spec ( {\bf A} (G))} (d(d \lambda {}^2 +4-d)-4 \lambda {}_A \lambda ) \\ 
\  &   &                \\ 
\  & = &  ( \lambda {}^2 -( \frac{4-d}{d} )^2 )^{m-n} / d^{2n} 
\prod_{ \lambda {}_A \in Spec ( {\bf A} (G))} ( d^2 \lambda {}^2 -4 \lambda {}_A \lambda +d(4-d)) . 
\end{array}
\]
Thus, solving 
\[
 d^2 \lambda {}^2 -4 \lambda {}_A \lambda +d(4-d)=0 ,
\]
we obtain  
\[
\lambda = \frac{2 \lambda {}_A \pm \sqrt{4 \lambda {}_A {}^2 - d^3 (4-d)}}{ d^2 } .
\]
$\Box$

In the case of $d=4$, we consider ${\bf P} =( P_{ef} )_{e,f \in D(G)} $ be the transition probability matrix of the CRW 
with respect to the Grover matrixon a $d$-regular graph $G$.
If $t(f)=o(e)$ and $f \not= e^{-1}$, then $ P_{ef} =4/d^2=4/4^2=1/4$. 
If $f = e^{-1} $, then $ P_{ef} =4/d^2 - (4/d - 1) = 4/4^2 - (4/4 - 1) = 1/4$. 
Thus, this CRW is considered to be a simple random walk on $G$ which the particle moves over each arc in terms of the same probability. 
Furthermore, an $n \times n$ Hadamard matrix is a unitary matrix whose elements have the absolute value $1/ \sqrt{n} $(see [2]). 
The Grover matrix of a $d$-regular graph is an Hadamard matrix if and only if $d=4$.   

\section{An application to the correlated random walk on a semiregular bipartite graph} 

We present spectra for the transition probability matrix of the correlated random walk on a 
semiregular bipartite graph. 
Hashimoto [4] presented a determinant expression for the Ihara zeta function 
of a semiregular bipartite graph. 
We use an analogue of the method in the proof of Hashimoto's result.

A bipartite graph $G=( V_1 ,V_2 )$ is called {\em $( q_1 , q_2 )$-semiregular} if 
$ \deg {}_G v = q_i $ for each $v \in V_i (i=1,2)$. 
For a $( q_1 +1, q_2 +1 )$-semiregular bipartite graph $G=( V_1 ,V_2 )$, 
let $G^{[i]} $ be the graph with vertex set $V_i $ and edge set 
$ \{ P : \ {\rm reduced \  path} \mid \  
\mid P \mid =2; o(P), t(P) \in V_i \}$ for $i=1,2$. 
Then $G^{[1]} $ is $( q_1 +1) q_2 $-regular, and $G^{[2]} $ 
is $( q_2 +1) q_1 $-regular.

\begin{theorem}
Let $G=(V, W)$ be a connected $(r,s)$-semiregular bipartite 
graph with $ \nu $ vertices and $ \epsilon $ edges. 
Set $ \mid V \mid =m$ and $ \mid W \mid =n(m \leq n)$.
Furthermore, let ${\bf P} $ be the transition probability matrix of the CRW with respect to the Grover matrix of $G$, and 
\[
Spec ( {\bf A} (G)= \{ \pm \lambda {}_1 , \cdots , \pm \lambda {}_{m} , 0, \ldots , 0 \} . 
\]  
Then 
\[
\det ( {\bf I}_{2 \epsilon } -u {\bf P} )=(1- u^2 (4/r-1)(4/s-1))^{\epsilon - \nu } (1- u^2 (4/r-1))^{n-m} 
\]
\[
\times \prod^m_{j=1} ((1- u^2 (4/s-1))(1- u^2 (4/r-1)) -16 \frac{\lambda {}^2_j }{r^2 s^2} u^2 ) .
\]
\end{theorem}

{\bf Proof}.  Let $e \in D(G)$. 
If $o(e) \in V$, then 
\[
d_{o(e)} =r, \  d_{t(e)} =s . 
\] 
Thus, we have 
\[
1-u^2 ( \frac{4}{ d_{o(e)}} -1)( \frac{4}{ d_{t(e)}} -1)= \frac{rs- u^2 (4-r)(4-s)}{rs} , 
\]
\[
\begin{array}{rcl} 
\ &  & ( {\bf A}_{CRW} )_{xy} = \frac{4/d^2_x }{1- u^2 (4/d_x -1)(4/d_y -1)} \\ 
\  &   &                \\ 
\  & = &  \left\{
\begin{array}{ll}
\frac{4s}{rs- u^2 (4-r)(4-s)} \frac{1}{r} & \mbox{if $(x,y) \in D(G)$ and $x \in V$, } \\
\frac{4r}{rs- u^2 (4-r)(4-s)} \frac{1}{s} & \mbox{if $(x,y) \in D(G)$ and $x \in W$, } 
\end{array}
\right. 
\end{array}
\]
and 
\[ 
\begin{array}{rcl} 
\ &  & ( {\bf D}_{CRW} )_{xx}= \sum_{o(e)=x} \frac{4/d^2_x (4/d_{t(e)} -1)}{1- u^2 (4/d_x -1)(4/d_{t(e)} -1)} \\ 
\  &   &                \\ 
\  & = &  \left\{
\begin{array}{ll}
r \cdot \frac{4(4-s)}{r(rs- u^2 (4-r)(4-s))} = \frac{4(4-s)}{rs- u^2 (4-r)(4-s))} & \mbox{if $x \in V$, } \\  
s \cdot \frac{4(4-r)}{s(rs- u^2 (4-r)(4-s))} = \frac{4(4-r)}{rs- u^2 (4-r)(4-s))} & \mbox{if $x \in W$. }
\end{array}
\right. 
\end{array}
\]

Next, let $V= \{ v_1, \cdots , v_m \} $ and 
$W= \{ w_1, \cdots , w_n \} $. 
Arrange vertices of $G$ as follows:
$v_1, \cdots , v_m ; w_1, \cdots , w_n $. 
We consider the matrix ${\bf A} = {\bf A} (G)$ under this order. 
Then, let 
\[
{\bf A} = 
\left[
\begin{array}{cc}
{\bf 0} & {\bf E} \\
{}^t {\bf E} & {\bf 0} 
\end{array}
\right]
.  
\]

Since $ {\bf A} $ is symmetric, there exists an orthogonal matrix 
${\bf F} \in O(n)$ such that 
\[
{\bf E} {\bf F} = 
\left[
\begin{array}{cc}
{\bf R} & {\bf 0} 
\end{array}
\right]
=
\left[
\begin{array}{cccccc}
\mu {}_1 &  & 0 & 0 & \cdots & 0 \\
 & \ddots &  & \vdots &  & \vdots \\
\star &  & \mu {}_{m} & 0 & \cdots & 0 
\end{array}
\right]
.
\]

Now, let 
\[
{\bf H} = 
\left[
\begin{array}{cc}
{\bf I}_{m} & {\bf 0} \\
{\bf 0} & {\bf F}
\end{array}
\right]
.
\]
Then we have 
\[
{}^t {\bf H} {\bf A} {\bf H} = 
\left[
\begin{array}{ccc}
{\bf 0} & {\bf R} & {\bf 0} \\
{}^t {\bf R} & {\bf 0} & {\bf 0} \\
{\bf 0} & {\bf 0} & {\bf 0} 
\end{array}
\right]
.
\]
Furthermore, let 
\[
\alpha =4/(rs- u^2 (4-r)(4-s)) . 
\]
Then we have 
\[
{\bf A}_{CRW} = 
\left[
\begin{array}{cc}
{\bf 0} & \alpha s/r {\bf E} \\
\alpha r/s \ {}^t {\bf E} & {\bf 0} 
\end{array}
\right]
, 
\]
and 
\[
{\bf D}_{CRW} = 
\left[
\begin{array}{cc}
\alpha (4-s) {\bf I}_m & {\bf 0} \\
{\bf 0} & \alpha (4-r) {\bf I}_n  
\end{array}
\right]
. 
\]
Thus, we have  
\[
{}^t {\bf H} {\bf A}_{CRW} {\bf H} =  
\left[
\begin{array}{ccc}
{\bf 0} & \alpha s/r {\bf R} & {\bf 0} \\
\alpha r/s \ {}^t {\bf R} & {\bf 0} & {\bf 0} \\
{\bf 0} & {\bf 0} & {\bf 0} 
\end{array}
\right]
\]
and 
\[
{}^t {\bf H} {\bf D}_{CRW} {\bf H} = 
\left[
\begin{array}{cc}
\alpha (4-s) {\bf I}_m & {\bf 0} \\
{\bf 0} & \alpha (4-r) {\bf I}_n  
\end{array}
\right]
. 
\]

By Theorem 4,  
\[
\begin{array}{rcl}
\  &   & \det ( {\bf I}_{2 \epsilon } -u {\bf P} )= \frac{(rs- u^2 (4-r)(4-s))^{\epsilon} }{ r^{\epsilon} s^{\epsilon} } 
\det ( {\bf I}_{\nu } -u {\bf A}_{CRW} + u^2 {\bf D}_{CRW} ) \\ 
\  &   &                \\ 
\  & = & \frac{(rs- u^2 (4-r)(4-s))^{\epsilon} }{ r^{\epsilon} s^{\epsilon} }  
\det 
\left( \left[
\begin{array}{ccc}
{\bf I}_m + \alpha (4-s)u^2 {\bf I}_m & -\alpha su/r {\bf R} & {\bf 0} \\
- \alpha ru/s\ {}^t {\bf R} & {\bf I}_m + \alpha (4-r)u^2 {\bf I}_m & {\bf 0} \\ 
{\bf 0} & {\bf 0} & {\bf I}_{n-m} + \alpha (4-r)u^2 {\bf I}_{n-m} 
\end{array}
\right] \right) 
\\
\  &   &                \\ 
\  & = & \frac{(rs- u^2 (4-r)(4-s))^{\epsilon} }{ r^{\epsilon} s^{\epsilon} } 
(1+ \alpha (4-r) u^2 )^{n-m} \\ 
\  &   &                \\ 
\  & \times & 
\det 
\left( \left[
\begin{array}{cc}
(1+ \alpha (4-s) u^2 ) {\bf I}_{m} & - \alpha su/r {\bf R} \\
- \alpha ru/s \  {}^t {\bf R} & (1+ \alpha (4-r) u^2 ) {\bf I}_{m} 
\end{array}
\right] \right) 
\cdot 
\det 
\left( \left[
\begin{array}{cc}
{\bf I}_{m} & \frac{1}{1+ \alpha (4-s) u^2 } \frac{\alpha su}{r} {\bf R} \\
{\bf 0} & {\bf I}_{m} 
\end{array}
\right] \right) 
\\
\  &   &                \\ 
\  & = & \frac{(rs- u^2 (4-r)(4-s))^{\epsilon +m-n} }{ r^{\epsilon} s^{\epsilon} } 
(rs+ u^2 (4-r)s)^{n-m} \\ 
\  &   &                \\ 
\  & \times & 
\det 
\left( \left[
\begin{array}{cc}
(1+ \alpha (4-s) u^2 ) {\bf I}_{m} & {\bf 0} \\
- \alpha ru/s \ {}^t {\bf R} & (1+ \alpha (4-r) u^2 ) {\bf I}_{m} 
- \frac{ \alpha {}^2 u^2 }{1+\alpha (4-s) u^2} \ {}^t {\bf R} {\bf R}   
\end{array}
\right] \right) 
\\
\  &   &                \\ 
\  & = & \frac{(rs- u^2 (4-r)(4-s))^{\epsilon +m-n} }{ r^{\epsilon} s^{\epsilon} } 
(rs+ u^2 (4-r)s)^{n-m} \\ 
\  &   &                \\ 
\  & \times & (1+ \alpha (4-s) u^2 )^m  
\det ((1+ \alpha (4-r) u^2 ) {\bf I}_{m} -\frac{\alpha {}^2 u^2 }{1+\alpha (4-s) u^2} \ {}^t {\bf R} {\bf R} ) \\ 
\  &   &                \\ 
\  & = &  \frac{(rs- u^2 (4-r)(4-s))^{\epsilon +m-n} }{ r^{\epsilon} s^{\epsilon} } 
(rs+ u^2 (4-r)s)^{n-m} \\ 
\  &   &                \\ 
\  & \times &   
\det ((1+ \alpha (4-s) u^2 )(1+ \alpha (4-r) u^2 ) {\bf I}_{m} - \alpha {}^2 u^2 \ {}^t {\bf R} {\bf R} )
.  
\end{array}
\]

Since ${\bf A} $ is symmetric, 
${}^t {\bf R} {\bf R} $ is symmetric and 
positive semi-definite, i.e., 
the eigenvalues of ${}^t {\bf R} {\bf R} $ are 
of form: 
\[
\lambda {}^2_1 , \cdots , \lambda {}^2_{m} 
( \lambda {}_1 , \cdots , \lambda {}_{m} \geq 0). 
\]
Furthermore, we have 
\[
\det ( \lambda {\bf I}_{\nu } - {\bf A} (G))= \lambda {}^{n-m} \det ( \lambda {}^2 - {}^t {\bf R} {\bf R} ) , 
\]
and so, 
\[
Spec ( {\bf A} (G)= \{ \pm \lambda {}_1 , \cdots , \pm \lambda {}_{m} , 0, \ldots , 0 \} . 
\]
Therefore it follows that 
\[
\begin{array}{rcl}
\  &   & \det ( {\bf I}_{2 \epsilon } -u {\bf P} ) \\ 
\  &   &                \\ 
\  & = & \frac{(rs- u^2 (4-r)(4-s))^{\epsilon +m-n} }{ r^{\epsilon} s^{\epsilon} } 
(rs+ u^2 (4-r)s)^{n-m} \\ 
\  &   &                \\ 
\  & \times &   
\prod^m_{j=1} ((1+ \alpha (4-s) u^2 )(1+ \alpha (4-r) u^2 ) {\bf I}_{m} - \alpha {}^2 \lambda {}^2_j u^2 ) \\ 
\  &   &                \\ 
\  & = & \frac{(rs- u^2 (4-r)(4-s))^{\epsilon +m-n} }{ r^{\epsilon} s^{\epsilon} } 
(rs+ u^2 (4-r)s)^{n-m} \\ 
\  &   &                \\ 
\  & \times &   
\prod^m_{j=1} ( \frac{rs+ u^2 (4-s)r}{rs- u^2 (4-r)(4-s)} \frac{rs+ u^2 (4-r)s}{rs- u^2 (4-r)(4-s)} \\ 
\  &   &                \\ 
\  & - & \lambda {}^2_j \frac{16 u^2 }{(rs- u^2 (4-r)(4-s))^2 } ) \\
\  &   &                \\ 
\  & = & \frac{(rs- u^2 (4-r)(4-s))^{\epsilon -m-n} }{ r^{\epsilon} s^{\epsilon} }  
(rs+ u^2 (4-r)s)^{n-m} \\ 
\  &   &                \\ 
\  & \times &   
\prod^m_{j=1} (rs(s+ u^2 (4-s))(r+ u^2 (4-r)) -16 \lambda {}^2_j u^2 ) \\
\  &   &                \\ 
\  & = & (1- u^2 (4/r-1)(4/s-1))^{\epsilon - \nu } (1+ u^2 (4/r-1))^{n-m} \\ 
\  &   &                \\ 
\  & \times &   
\prod^m_{j=1} ((1+ u^2 (4/s-1))(1+ u^2 (4/r-1)) -16 \frac{\lambda {}^2_j }{r^2 s^2} u^2 ) .   
\end{array}
\]
$\Box$

Now, let $u=1/ \lambda $. 
Then we obtain the following result.

\begin{corollary}
Let $G=(V, W)$ be a connected $(r,s)$-semiregular bipartite 
graph with $ \nu $ vertices and $ \epsilon $ edges. 
Set $ \mid V \mid =m$ and $ \mid W \mid =n(m \leq n)$.
Furthermore, let ${\bf P} $ be the transition probability matrix of the CRW with respect to the Grover matrix and 
\[
Spec ( {\bf A} (G)= \{ \pm \lambda {}_1 , \cdots , \pm \lambda {}_{m} , 0, \ldots , 0 \} . 
\]  
Then 
\[
\det ( \lambda {\bf I}_{2 \epsilon } - {\bf P} )=( \lambda {}^2 -(4/r-1)(4/s-1))^{\epsilon - \nu } ( \lambda {}^2 +(4/r-1))^{n-m} 
\]
\[
\times \prod^m_{j=1} (( \lambda {}^2 +(4/s-1))( \lambda {}^2 +(4/r-1)) -16 \frac{\lambda {}^2_j }{r^2 s^2} \lambda {}^2 ) .
\]
\end{corollary}

By Corollary 4,  we obtain the spectra for the transition probability matrix ${\bf P} $ of the CRW with respect to the Grover matrix 
of a semiregular bipartite graph.

\begin{corollary}
Let $G=(V, W)$ be a connected $(r,s)$-semiregular bipartite 
graph with $ \nu $ vertices and $ \epsilon $ edges. 
Set $ \mid V \mid =m$ and $ \mid W \mid =n(m \leq n)$.
Furthermore, let ${\bf P} $ be the transition probability matrix of the CRW with respect to the Grover matrix and 
\[
Spec ( {\bf A} (G)= \{ \pm \lambda {}_1 , \cdots , \pm \lambda {}_{m} , 0, \ldots , 0 \} . 
\] 
Then the transition matrix ${\bf P} $ has $2 \epsilon $ eigenvalues of the form 
\begin{enumerate} 
\item $4m$ eigenvalues: 
\[  
\lambda = \pm \sqrt{ \frac{2r^2 s^2 - 4rs^2-4r^2 s+16 \lambda {}^2_j 
\pm \sqrt{ (2r^2 s^2 - 4rs^2-4r^2 s+16 \lambda {}^2_j )^2 -4 r^3 s^3 (4-r)(4-s)} }{2 r^2 s^2 } } ;   
\] 
\item $2n-2m$ eigenvalues: 
\[
\lambda = \pm i \sqrt{ \frac{4}{r} -1} ; 
\]
\item $2( \epsilon - \nu)$ eigenvalues: 
\[
\lambda = \pm \sqrt{( \frac{4}{r} -1)( \frac{4}{s} -1)} . 
\] 
\end{enumerate} 
\end{corollary}

{\bf Proof}.  Solving 
\[
( \lambda {}^2 +(4/s-1))( \lambda {}^2 +(4/r-1)) -16 \frac{\lambda {}^2_j }{r^2 s^2} \lambda {}^2  =0 , 
\]
i.e., 
\[
\lambda {}^4 +( \frac{4}{r} + \frac{4}{s} -2- \frac{16 \lambda {}^2_j }{ r^2 s^2 } ) \lambda {}^2 +( \frac{4}{r} -1)( \frac{4}{s} -1)=0 , 
\]
we obtain 
\[ 
\lambda = \pm \sqrt{ \frac{1}{2} ((2- \frac{4}{r} - \frac{4}{s} + \frac{16 \lambda {}^2_j }{ r^2 s^2 } ) \pm 
\sqrt{ (2- \frac{4}{r} - \frac{4}{s} + \frac{16 \lambda {}^2_j }{ r^2 s^2 } )^2 -4 ( \frac{4}{r} -1)( \frac{4}{s} -1) } } , 
\]
i.e., 
\[  
\lambda = \pm \sqrt{ \frac{2r^2 s^2 - 4rs^2-4r^2 s+16 \lambda {}^2_j 
\pm \sqrt{ (2r^2 s^2 - 4rs^2-4r^2 s+16 \lambda {}^2_j )^2 -4 r^3 s^3 (4-r)(4-s)} }{2 r^2 s^2 } } ;   
\] 
$\Box$

\section{Another type of the correlated random walk on a cycle graph}

The CRW is defined by the following transition probability matrix ${\bf P} $ 
on the one-dimensional lattice: 
\[
{\bf P} =  
\left[
\begin{array}{cc}
a & b \\
c & d 
\end{array}
\right]
, 
\]
where 
\[
a+c=b+d=1, \ a,b,c,d \in [0,1] . 
\] 
As for the CRW, see [7,9], for example. 

We formulate a CRW on the arc set of a graph with respect to the above matrix ${\bf P} $. 
The cycle graph is a connected 2-regular graph. 
Let $C_n $ be the cycle graph with $n$ vertices and $n$ edges. 
Furthermore, let $V(C_n )= \{ v_1 , \ldots , v_n \} $ and $e_j =(v_j , v_{j+1} ) (1 \leq j \leq n)$, 
where the subscripts are considered by modulo $m$.  
Then we introduce a $2n \times 2n$ matrix ${\bf U} =( U_{ef} )_{e,f \in D(C_n )} $ as follows: 
\[ 
U_{ef} = \left\{
\begin{array}{ll}
d & \mbox{if $t(f)=o(e)$, $f \neq e^{-1}$ and $f=e_j $, } \\
b & \mbox{if $f= e^{-1} $ and $f=e_{j} $, } \\
a & \mbox{if $t(f)=o(e)$, $f \neq e^{-1} $ and $f=e^{-1}_{j} $, } \\
c & \mbox{if $f= e^{-1} $ and $f=e^{-1}_j $, } \\
0 & \mbox{otherwise.}
\end{array}
\right. 
\] 
Note that 
${\bf U} $ is be able to write as follows: 
\[
{\bf U} =  
\left[
\begin{array}{cc}
d {\bf Q}^{-1} & c {\bf I}_n \\
b {\bf I}_n & a {\bf Q}  
\end{array}
\right]
, 
\]
where ${\bf Q} = {\bf P}_{\sigma } $ is the permutation matrix of $\sigma =(12 \ldots n)$.  
The CRW with ${\bf U} $ with a transition probability matrix is called the {\em second type of CRW} on $C_n $ 
with respect to the above matrix ${\bf P} $. 

Now, we define a function $w: D(C_n) \longrightarrow \mathbb{R} $ as follows: 
\[
w(e)= \left\{
\begin{array}{ll}
d & \mbox{if $e=e_j \ (1 \leq j \leq n)$, } \\ 
a & \mbox{if $e= e^{-1}_j \ (1 \leq j \leq n)$. } 
\end{array}
\right. 
\]
Furthermore, let an $n \times n$ matrix $ {\bf W} (C_n )=( w_{uv} )_{u,v \in V(C_n )} $ as folloows: 
\[
w_{uv} = \left\{
\begin{array}{ll}
w(u,v) & \mbox{if $(u,v) \in D(C_n)$, } \\ 
0 & \mbox{otherwise. }
\end{array}
\right. 
\]

The characteristic polynomial of ${\bf U} $ is given as follows.

\begin{theorem}
Let $C_n $ be the cycle graph with $n$ vertices, and ${\bf U} $ the transition probability matrix of 
the second type of CRW on $C_n $. 
Then 
\[
\det ( \lambda {\bf I}_{2n} - {\bf U} )= \det (( \lambda {}^2 +(ad-bc)) {\bf I}_n - \lambda {\bf W} (C_n )) . 
\]
\end{theorem}

{\bf Proof}.  At first, we consider two $2n \times 2n$ matrices $2n \times 2n$ matrices ${\bf B} =( B_{ef} )_{e,f \in D(C_n)} $ and 
${\bf J} =( J_{ef} )_{e,f \in D(C_n)} $ as follows: 
\[
B_{ef} = \left\{
\begin{array}{ll}
w(f) & \mbox{if $t(e)=o(f)$, } \\ 
0 & \mbox{otherwise, }
\end{array}
\right. 
\ 
J_{ef} = \left\{
\begin{array}{ll}
b-a & \mbox{if $f=e^{-1}$ and $e=e_j$, } \\ 
c-d & \mbox{if $f=e^{-1}$ and $e=e^{-1}_j$, } \\ 
0 & \mbox{otherwise. }
\end{array}
\right.  
\]
Then we have 
\[
{\bf U} = {}^t {\bf B} + {}^t {\bf J} .  
\]

Now, we define two $2n \times n$ matrices ${\bf K} =( K_{ev} )_{e \in D(C_n) ; v \in V(C_n )} $ and 
${\bf L} =( L_{ev} )_{e \in D(C_n) ; v \in V(C_n )} $ as follows: 
\[
K_{ev} = \left\{
\begin{array}{ll}
1 & \mbox{if $t(e)=v$, } \\ 
0 & \mbox{otherwise, }
\end{array}
\right. 
\ 
L_{ev} = \left\{
\begin{array}{ll}
w(e) & \mbox{if $o(e)=v$, } \\ 
0 & \mbox{otherwise. }
\end{array}
\right.
\]
Then we have 
\[
{\bf K} \ {}^t {\bf L} = {\bf B} , \ {}^t {\bf L} {\bf K} = {\bf W} (C_n ) . 
\]
If ${\bf A} $ and ${\bf B} $ are an $m \times n$ matrix and an $n \times m$ matrix, respectively, then 
we have 
\[
\det ( {\bf I}_m - {\bf A} {\bf B} ) = \det ( {\bf I}_n - {\bf B} {\bf A} ) . 
\]
Thus, 
\[
\begin{array}{rcl} 
\det ( {\bf I}_{2n} -u {\bf U} ) & = & \det ( {\bf I}_{2n} -u( {}^t {\bf B} + {}^t {\bf J} )) \\ 
\  &   &                \\ 
\  & = & \det ( {\bf I}_{2n} -u( {\bf B} + {\bf J} )) \\ 
\  &   &                \\ 
\  & = & \det ( {\bf I}_{2n} -u {\bf J} -u {\bf B} ) \\ 
\  &   &                \\ 
\  & = & \det ( {\bf I}_{2n} -u {\bf J} -u {\bf K} \ {}^t {\bf L} ) \\ 
\  &   &                \\ 
\  & = & \det ( {\bf I}_{2n} -u {\bf K} \ {}^t {\bf L} ( {\bf I}_{2n} -u {\bf J} )^{-1} ) 
\det ( {\bf I}_{2n} -u {\bf J} ) \\ 
\  &   &                \\
\  & = & \det ( {\bf I}_{n} -u \ {}^t {\bf L} ( {\bf I}_{2n} -u {\bf J} )^{-1} {\bf K} ) \det ( {\bf I}_{2n} -u {\bf J}) .  
\end{array}
\]

But, we have 
\[
\begin{array}{rcl} 
\  &  & \det ( {\bf I}_{2n} -u {\bf J} ) \\ 
\  &   &                \\ 
\  & = &  
\left[
\begin{array}{cc}
{\bf I}_n & -(b-a)u {\bf I}_n \\
-(c-d)u {\bf I}_n & {\bf I}_n  
\end{array}
\right]
\cdot 
\left[
\begin{array}{cc}
{\bf I}_n & (b-a)u {\bf I}_n \\
{\bf 0} & {\bf I}_n  
\end{array}
\right]
\\
\  &   &                \\ 
\  & = & 
\left[
\begin{array}{cc}
{\bf I}_n & {\bf 0} \\
-(c-d)u {\bf I}_n & {\bf I}_n -u^2 (b-a)(c-d) {\bf I}_n  
\end{array}
\right]
\\ 
\  &   &                \\ 
\  & = & (1-(a-c)(d-b) u^2 )^n . 
\end{array}
\]
Furthermore, we have 
\[
( {\bf I}_{2n} -u {\bf J} )^{-1} = \frac{1}{1-(a-b)(d-c) u^2 } ( {\bf I}_{2n} +u {\bf J} ) . 
\]
Therefore, it follows that 
\[
\begin{array}{rcl} 
\  &  & \det ( {\bf I}_{2n} -u {\bf U} ) \\ 
\  &   &                \\ 
\  & = & (1-(a-b)(d-c) u^2 )^n \det ( {\bf I}_{n} -u/(1-(a-b)(d-c) u^2 ) \ {}^t {\bf L} ( {\bf I}_{2n} +u {\bf J} ) {\bf K} ) \\ 
\  &   &                \\ 
\  & = & \det ((1-(a-b)(d-c) u^2 ) {\bf I}_{n} -u \ {}^t {\bf L} {\bf K} -u \ {}^t {\bf L} {\bf J} {\bf K} ) \\ 
\  &   &                \\ 
\  & = & \det ((1-(a-b)(d-c) u^2 ) {\bf I}_{n} -u {\bf W} (C_n )-u^2 \ {}^t {\bf L} {\bf J} {\bf K} ) . 
\end{array}
\]

The matrix ${}^t {\bf L} {\bf J} {\bf K}$ is a diagonal, and its $(v_i , v_i)$ entry is equal to  
\[
(c-d) w(e^{-1}_{i-1} )+(b-a) w( e_i )=(c-d)a+(b-a)d=c+bd-2ad . 
\]
That is, 
\[
{}^t {\bf L} {\bf J} {\bf K} =(ab+cd-2ad) {\bf I}_n . 
\]
Thus, 
\[
\begin{array}{rcl} 
\  &  & \det ( {\bf I}_{2n} -u {\bf U} ) \\ 
\  &   &                \\ 
\  & = &  \det ((1-(a-b)(d-c) u^2 ) {\bf I}_{n} -u {\bf W} (C_n )-u^2 (ac+bd-2ad) {\bf I}_n ) \\ 
\  &   &                \\ 
\  & = &  \det (((1+(ad-bc) u^2 ) {\bf I}_{n} -u {\bf W} (C_n )) . 
\end{array}
\]

Substituting $u=1/ \lambda$, the result follows. 
$\Box$

By Theorem 7, we obtain the spectra for the transition probability matrix ${\bf U} $ of 
the second type of the CRW on $C_n $.
The matrix ${\bf W} (C_n )$ is given as follows: 
\[ 
{\bf W} (C_n )=  
\left[
\begin{array}{cccccc}
0 & d & 0 & \ldots &  & a \\
a & 0 & d & \ldots &  & 0 \\
\vdots & \vdots &  & \ddots &  & \\
0 & 0 & 0 & \ldots & 0 & d \\
d & 0 & 0 & \ldots & a & 0 
\end{array}
\right]
, 
\]

\begin{corollary}
Let $C_n $ be the cycle graph with $n$ vertices, and ${\bf U} $ the transition probability matrix of 
the second type of CRW on $C_n $. 
Then the transition probability matrix ${\bf U} $ has $2n$ eigenvalues of the form 
\[
\lambda = \frac{\mu \pm \sqrt{ \mu {}^2 -4(ad-bc)}}{2} , \ \mu \in Spec ( {\bf W} (C_n )) .  
\]
\end{corollary}

{\bf Proof}.  At first, we have 
\[
\det ( {\bf I}_{2n} -u {\bf U} )= \prod_{ \mu \in  Spec ( {\bf W} (C_n ))} ( \lambda {}^2 - \mu \lambda +(ad-bc)) . 
\]

Solving 
\[
\lambda {}^2 - \mu \lambda +(ad-bc)=0 , 
\]
we obtain 
\[ 
\lambda = \frac{\mu \pm \sqrt{ \mu {}^2 -4(ad-bc)}}{2} . 
\] 
$\Box$

Now, we consider the case of $a=b=c=d=1/2$. 
Then the matrix ${\bf W} (C_n )$ is equal to 
\[
{\bf W} (C_n )= \frac{1}{2} {\bf A} (C_n ) . 
\]

By Corollary 6, we obtain the spectra for the transition probability matrix ${\bf U} $ of the second type of 
CRW on $C_n $.

\begin{corollary}
Let $C_n $ be the cycle graph with $n$ vertices, and ${\bf U} $ the transition probability matrix of 
the second type of the CRW on $C_n $. 
Assume that $a=b=c=d=1/2$.
Then the transition probability matrix ${\bf U} $ has $n$ eigenvalues of the form 
\[
\lambda = \cos \theta {}_j , \ \theta {}_j = \frac{2 \pi j}{n} (j=0,1, \ldots , n-1) \  \  \  (*).  
\]
The remaining $n$ eigenvalues of ${\bf U} $ are $0$ with multiplicities $n$. 
\end{corollary}

{\bf Proof}. It is known that the spectrum of ${\bf A} (C_n )$ are 
\[ 
2 \cos \theta {}_j , \ \theta {}_j = \frac{2 \pi j}{n} (j=0,1, \ldots , n-1) .
\]
$\Box$

Note that the spectrum of (*) are those of the transition probability matrix of the simple random walk 
on a cycle graph $C_n $.

We can generalize the result for $a=b=c=d=1/2$ on $C_n $ to a $d$-regular graph$(d \geq 2)$.  
Let $G$ be a connected $d$-regular graph with $n$ vertices and $m$ edges.  
Furthermore, let ${\bf P} $ be the $d \times d$ matrix as follows: 
\[
{\bf P} = \frac{1}{d} {\bf J}_d , 
\]
where ${\bf J}_d $ is the matrix whose elements are all one. 
Let ${\bf U} =( U_{ef} )_{e,f \in D(G)}$ be the the transition probability matrix of a CRW on $G$ 
with respect to ${\bf P} $. 
Then we have 
\[
U_{ef} = \left\{
\begin{array}{ll}
1/d & \mbox{if $t(e)=o(f)$, } \\ 
0 & \mbox{otherwise, }
\end{array}
\right. 
\] 
and so, 
\[
{\bf U} = \frac{1}{d} {\bf B} . 
\]

Similarly to The proof of Theorem 7, we obtain the following result.

\begin{theorem}
Let $G$ be a connected $d$-regular graph with $n$ vertices and $m$ edges.  
Furthermore, let ${\bf U} $ the transition probability matrix of 
the CRW on $G$ with respect to ${\bf P} =1/d {\bf J}_d $. 
Then 
\[
\det ( \lambda {\bf I}_{2m} - {\bf U} )= \lambda {}^{2m-n} \det ( \lambda {\bf I}_n -\frac{1}{d} {\bf A} (G)) . 
\]
\end{theorem}

Thus,

\begin{corollary}
Let $G$ be a connected $d$-regular graph with $n$ vertices and $m$ edges.  
Furthermore, let ${\bf U} $ the transition probability matrix of 
the CRW on $G$ with respect to ${\bf P} =1/d {\bf J}_d $. 
Then the transition probability matrix ${\bf U} $ has $n$ eigenvalues of the form 
\[
\lambda = \frac{1}{d} \lambda {}_A , \ \lambda {}_A \in Spec ( {\bf A} (G)) . 
\]
The remaining $2(m-n)$ eigenvalues of ${\bf U} $ are $0$ with multiplicities $2m-n$. 
\end{corollary}

\end{document}